\documentclass[12pt]{amsart}

\usepackage{times,epsf,amssymb,amsmath,hyperref, rlepsf, color}

\begin{document}

\newtheorem{thm}{Theorem}[section]
\newtheorem{lem}[thm]{Lemma}
\newtheorem{cor}[thm]{Corollary}
\newtheorem{conj}[thm]{Conjecture}
\newtheorem{question}[thm]{Question}

\theoremstyle{definition}
\newtheorem{defn}[thm]{\bf{Definition}}

\theoremstyle{remark}
\newtheorem{rmk}[thm]{Remark}

\def\square{\hfill${\vcenter{\vbox{\hrule height.4pt \hbox{\vrule width.4pt height7pt \kern7pt \vrule width.4pt} \hrule height.4pt}}}$}

\newenvironment{pf}{{\it Proof:}\quad}{\square \vskip 12pt}

\title[$H$-Surfaces in Homology Classes]{Constant Mean Curvature Surfaces in Homology Classes}

\author{Baris Coskunuzer}
\address{UT Dallas, Dept. Math. Sciences, Richardson, TX 75080}
\email{baris.coskunuzer@utdallas.edu}

\thanks{The author is partially supported by Simons Collaboration Grant, and Royal Society Newton Mobility Grant.}

\maketitle

%% User definitions:

\newcommand{\Si}{S^2_{\infty}({\Bbb H}^3)}
\newcommand{\SI}{S^n_{\infty}({\Bbb H}^{n+1})}
\newcommand{\PI}{\partial_{\infty}}

\newcommand{\BH}{\Bbb H}
\newcommand{\BHH}{{\Bbb H}^3}
\newcommand{\BR}{\Bbb R}
\newcommand{\BN}{\Bbb N}
\newcommand{\BC}{\Bbb C}
\newcommand{\BZ}{\Bbb Z}

\newcommand{\e}{\epsilon}

\newcommand{\wh}{\widehat}
\newcommand{\wt}{\widetilde}

\newcommand{\A}{\mathcal{A}}
\newcommand{\C}{\mathcal{C}}
\newcommand{\D}{\mathcal{D}}
\newcommand{\U}{\mathcal{U}}
\newcommand{\F}{\mathcal{F}}
\newcommand{\I}{\mathcal{I}}
\newcommand{\T}{\mathcal{T}}
\newcommand{\B}{\mathbf{B}}
\newcommand{\h}{\mathbf{h}}
\newcommand{\kk}{\mathbf{k}}

\begin{abstract}
We show the existence of constant mean curvature surfaces in the homology classes of closed $3$-manifolds.
\end{abstract}

\section{Introduction}

In this paper, we study the existence of constant mean curvature (CMC) surfaces in closed Riemannian $3$-manifolds. There are several results on the subject in various settings \cite{Hi, HK, DS, Gu, Ka, We,Ma,Me,MPT, Co1}. In particular, when $H=0$, this is a classical problem of the geometric analysis: Existence of minimal surfaces in $3$-manifolds.

One particular result in this context closely related to our discussion is the existence of area minimizing surfaces in a given homology class of a $3$-manifold. In particular, let $M$ be a closed Riemannian $3$-manifold. Let $S$ be an embedded surface in $M$ where $S$ has nontrivial homology, i.e. $[S]\neq 0\in H_2(M)$. Then, there exists a smoothly embedded area minimizing surface $\Sigma_0$ in $[S]$ by celebrated results of geometric measure theory \cite{Fe, HS}. 

In this paper, we will consider a generalization of this result to CMC surfaces. There are a couple of natural questions about this generalization.

\begin{question} Is there a CMC surface $\Sigma_H$ in $[S]$ for $H>0$?
	\end{question}

Then, the natural followup question:

\begin{question} For a given such $S\subset M$, for which $H$, there exists a CMC surface $\Sigma_H$ in $[S]$?
\end{question}

In this paper, we will address these two questions. Our main result is as follows:

\begin{thm}
	Let $M$ be a closed Riemannian $3$-manifold. Let $S$ be a closed embedded surface in $M$, where $S$ has nontrivial homology. Then, there exists $\widehat{H}_{[S]}\geq C_{[S]}\geq 0$ such that for any $H\in [0,\wh{H}_{[S]}]$, there exists a smoothly embedded $H$-surface $\Sigma_H$ in the homology class of $S$.
\end{thm}

Here, $C_{[S]}=\dfrac{g([S])-|\Sigma_0|}{||M||}$ where $|.|$ represents area, and $||.||$ represents volume. Also, we introduce a notion called {\em girth of a homology class}, $g([S])$. In particular,  $g([S])=\max\I_{[S]}$ where $\I_{[S]}$ is the isoperimetric profile of the homology class $[S]$, which is defined in the next section. 

Notice that our result is two-fold. First, we show that there is no gap in the CMC spectrum of the homology class ($[0,\wh{H}_{[S]}]$). The second is that we are giving a lower bound $C_{[S]}$ for $\wh{H}_{[S]}$ where $\wh{H}_{[S]}=\max\{H\mid \exists \Sigma_H\in [S]\}$. We will discuss how good is this lower bound in the final section.

One important remark at this point is the following. There are some trivial examples where there is no CMC surface in the homology class for $H\neq 0$. In particular, let $M=T\times S^1$ where $T$ is a closed surface. Let $M$ has the product metric. Then, by maximum principle, there is no CMC surface in the homology class of $T$ in $M$ for any $H>0$. On the other hand, for this example $C_{[T]}=0$ as the isoperimetric profile $\I_{[S]}$ is a constant function, and the girth of $[T]$, $g([T])=|T_0|$ in this particular example. So, even though our lower bound is not sharp, it gives a good lower bound in generic cases.

We will use two different techniques to get CMC surfaces in the proof of the main result. To show that there is no gap in the CMC spectrum of the homology class, we will use minimizing $H$-surfaces. On the other hand, in order to get a lower bound for $\wh{H}$, we will use isoperimetric surfaces. Meanwhile, we will introduce several notions like \textit{girth of a homology class} $g([S])$, \textit{isoperimetric profile of a homology class} $\I_{[S]}$. These are natural generalizations of their counterparts in the literature.

In the final section, we will discuss the correspondence between these minimizing $H$-surfaces and isoperimetric surfaces. In that section, we will also discuss further questions, and possible ways to improve this lower bound by using the recent techniques involving min-max surfaces.

%Organization of the paper is as follows: In the following section, we will give some definitions, and overview the related results which will be used in the following sections. In Section 3, we will prove our main result. In Section 4, we will give discuss further directions, and questions on the problem. In the final section, we have an appendix where we prove an interesting related result.

\subsection{Acknowledgements} I would like to thank Xin Zhou for very valuable conversations, and remarks.

\section{Preliminaries}

Let $M$ be a closed Riemannian $3$-manifold. Let $S$ be an embedded closed surface in $M$ with $[S]\neq 0\in H_2(M)$. Let $\Sigma_0$ be the area minimizing surface in $[S]$ \cite{Fe, HS}. 

Consider the noncomplete $3$-manifold $M-\Sigma_0$. Let $\wh{M}$ be the metric completion of $M-\Sigma_0$. In particular, $\wh{M}$ has two boundary components $\Sigma_0^+$ and $\Sigma_0^-$ which are exact copies of the area minimizing surface $\Sigma_0$. See Figure \ref{fig2} for an analogous picture in one dimension lower.

\subsection{Minimizing $H$-Surfaces} \label{minHsec} \

\begin{figure}[h]
	\begin{center}
		$\begin{array}{c@{\hspace{.4in}}c}

		\relabelbox  {\epsfxsize=1.7in \epsfbox{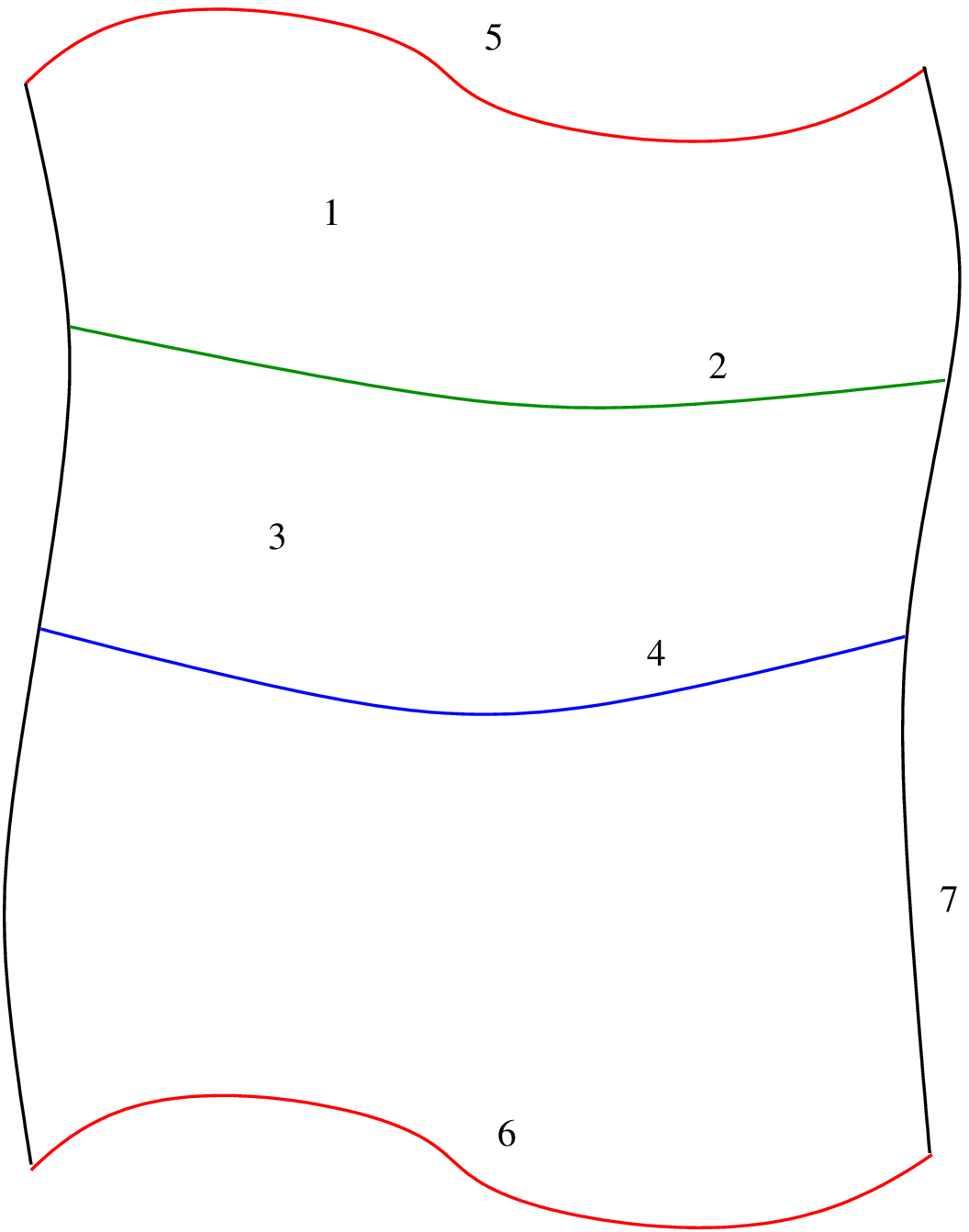}}
		\relabel{1}{$\Omega$} 
		\relabel{2}{\small $\Sigma$}  
		\relabel{3}{$\Omega'$} 
		\relabel{4}{\small $\wh{\Sigma}$} 
		\relabel{5}{\small $\Sigma_0^+$} 
		\relabel{6}{\small $\Sigma_0^-$} 	
		\relabel{7}{$\wh{M}$} 	   \endrelabelbox &
		
		\relabelbox  {\epsfxsize=1.7in \epsfbox{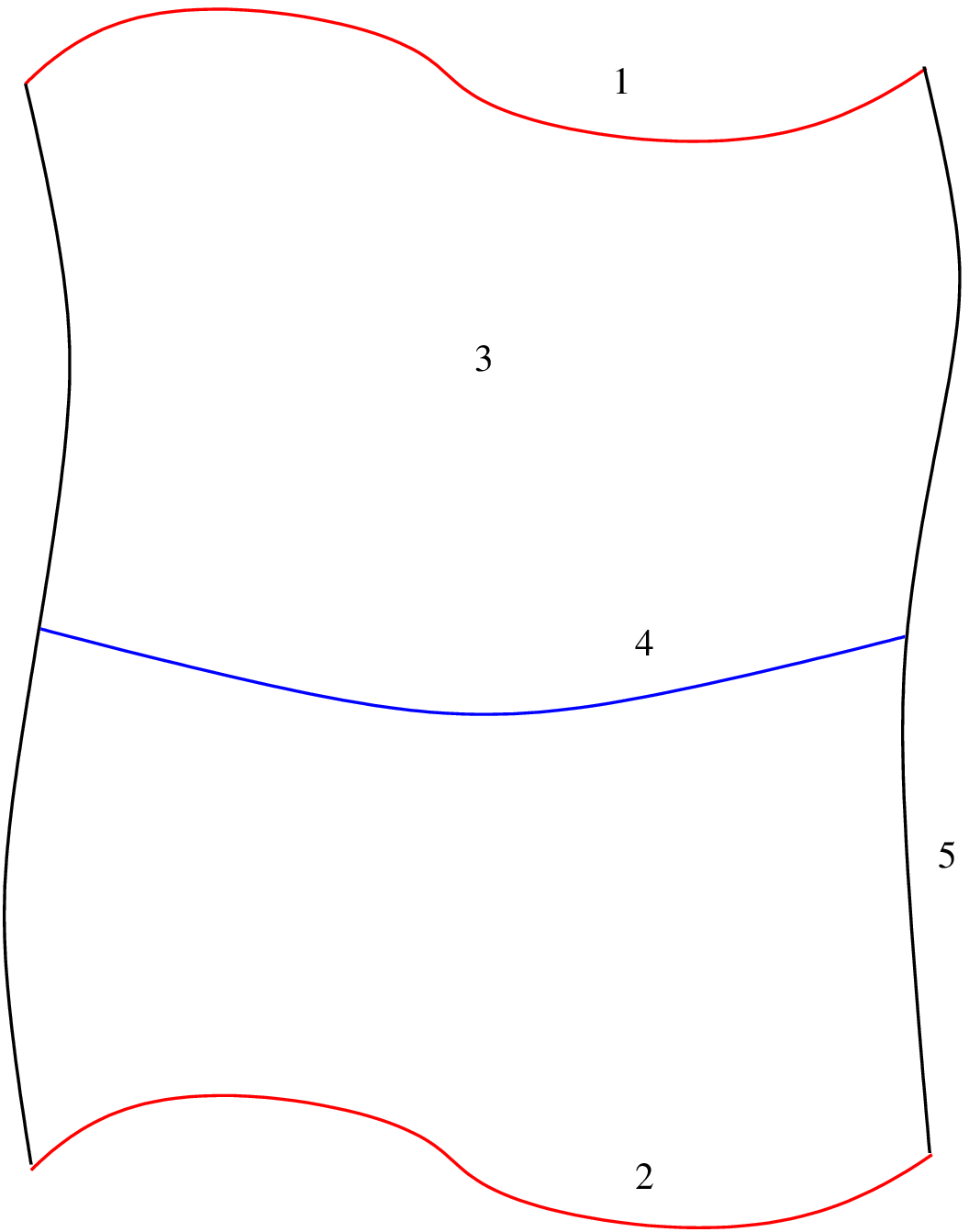}}
		\relabel{1}{\small $\Sigma_0^+$} 
		\relabel{2}{\small $\Sigma_0^-$}  
		\relabel{3}{$\Delta$} 
		\relabel{4}{\small $\T_K$}
		\relabel{5}{$\wh{M}$}	\endrelabelbox \\
		
	\end{array}$
	
\end{center}

\caption{ \label{fig1} \footnotesize In the figure left, we use $\wh{\Sigma}$ as a barrier surface to define $\A_H$. In the figure right, the least area separating surface in $\wh{M}$ with $\|\Delta\|=K$ gives the isoperimetric surface $\T_K\in [\Sigma_0]$.}
\end{figure}

Assuming the existence of a "nice" barrier surface $\wh{\Sigma}\in [S]$, we can define the following functional where the minimizers are CMC surfaces. Let $\wh{\Omega}$ be the region in $\wh{M}$ between $\wh{\Sigma}$ and $\Sigma_0^+$. Let $\Sigma$ be a surface in $\wh{\Omega}$ separating $\Sigma_0^+$ and $\wh{\Sigma}$. Let $\Omega'$ be the region between $\Sigma$ and $\wh{\Sigma}$ in $\wh{M}$. Let $\Omega$ be the region between $\Sigma$ and $\Sigma_0^+$ (See Figure \ref{fig1}-left). Notice that for any such $\Sigma$, $\Omega\cup\Omega'=\wh{\Omega}$. Now, define $$\A_H=|\Sigma|+2H||\Omega'||$$

This is an area functional with a volume constraint. The minimizers of $\A_H$ are $H$-surfaces. On the other hand, we can modify the functional as follows: 
$$\wh{\A}_H=|\Sigma|-2H||\Omega||$$

Notice that as $||\Omega||+||\Omega'||=||\wh{\Omega}||=C_0$, $\wh{\A}_H=\A_H-2HC_0$. Hence, $\A_H$ and $\wh{\A}_H$ differ by a constant. Therefore, minimizers of $\A_H$ and $\wh{\A}_H$ are same. We will use both functional interchangeably throughout the paper.

One final remark about the barrier surface $\wh{\Sigma}$ is the following. If $\wh{\Sigma}$ is a mean convex surface with $H(x)\geq H_0$ for any $x\in \wh{\Sigma}$, then for any $H\in [0,H_0)$ the minimizer $\Sigma_H$ would be disjoint from $\partial \wh{\Omega}=\Sigma_0^+\cup\wh{\Sigma}$ by the maximum principle. This would imply $\Sigma_H$ is a smoothly embedded surface by the regularity results of the geometric measure theory \cite{Fe}.

\subsection{Isoperimetric Surfaces} \

Let $M, S, \Sigma_0,\wh{M}$ be described as above. Define \textit{isoperimetric surfaces} $\T_K$ in $\wh{M}$ as follows. Let $\T$ be a surface in $\wh{M}$ separating $\Sigma_0^+$ and $\Sigma_0^-$. Let $\Delta$ be the region between $\T$ and $\Sigma_0^+$. 

Fix $K_0\in [0,||M||]$. Consider all such surfaces $\T$ with $||\Delta||=K_0$ (See Figure \ref{fig1}-right). Call the smallest area such surface $\T_{K_0}$ as the isoperimetric surface, i.e. $$|\T_{K_0}|=\inf_{||\Delta||=K_0} \{|\T|\}$$.

Now, we define the \textit{isoperimetric profile function} $\I_{[S]}$ of the homology class $[S]$ as follows.  Let $\T$ be an embedded surface in $\wh{M}$ separating $\Sigma_0^+$ and $\Sigma_0^-$. Let $\Delta$ be the region in $\wh{M}$ separated by $\Sigma_0^+$ and $\T$. Let $K\in [0,||M||]$. Then define 
$$\I_{[S]}(K)=inf\{|\T| \mid ||\Delta||=K \}$$

In \cite{Ro}, there are several results given for the isoperimetric profile of closed $3$-manifold $M$. By the proof of Lemma \ref{ctslem}, we can adapt these results to our setting $\wh{M}$. In particular, $\I_{[S]}(K)$ is a continuous function. Furthermore, for any $K\in [0,||M||]$, right and left derivatives of the function $\I_{[S]}(K)$ exists, i.e. $\I_{[S]}^{'+}(K)$ (right derivative) and $\I_{[S]}^{'-}(K)$ (left derivative). Moreover, there exists an isoperimetric surface $\T_K^\pm$ with mean curvature $H^\pm=\frac{\I_{[S]}^{'\pm}(K)}{2}$. By using this, we will get our lower bound $C_{[S]}$. Note that the derivative can be negative, which means the mean curvature vector changed direction.

\section{Existence of CMC Surfaces in Homology Classes}

In this section, we will prove our main result. In the first part, we will show that there is no gap in the CMC spectrum of the homology class. In the second part, we will give a lower bound for $\wh{H}_{[S]}$, the maximum $H$ for the CMC surfaces in the homology class $[S]$.

\subsection{No gap in the CMC Spectrum of a Homology Class.} \label{nogapsec} \

Let $M$ be a closed Riemannian $3$-manifold. Let $S$ be an embedded closed surface in $M$ with $[S]\neq 0\in H_2(M)$. Let $\Sigma_0$ be the area minimizing surface in $[S]$ \cite{Fe}. 

Consider the noncomplete $3$-manifold $M-\Sigma_0$. Let $\wh{M}$ be the metric completion of $M-\Sigma_0$. In particular, $\wh{M}$ has two boundary components $\Sigma_0^+$ and $\Sigma_0^-$ which are exact copies of the area minimizing surface $\Sigma_0$. 

In Figure \ref{fig2}, we gave an example of this in one dimension lower. In the figure, $M$ is a genus $3$ surface, and  $\Sigma_0$  was represented by an embedded geodesic curve (codimension-$1$) in its non-trivial homology class $[\Sigma_0]\neq 0\in H_1(M)$. Furthermore, $\wh{S}$ was represented by a homologous curve to $\Sigma_0$ in $M$. By cutting $M$ along $\Sigma_0$, we get a new manifold $\wh{M}$ which is a genus $2$ surface with two boundary components $\Sigma_0^+$ and $\Sigma_0^-$. Notice that $\wh{S}$ is a curve in $\wh{M}$ separating $\Sigma_0^+$ and $\Sigma_0^-$.

Notice that any surface $\Sigma$ in $\wh{M}$ separating $\Sigma_0^+$ and $\Sigma_0^-$ is a surface in the homology class of $[S]$ as $\Sigma\cup \Sigma_0^+$ separates a region $\Omega^+$ from $\wh{M}$. As $\Omega^+$ is also region in $M$ with $\partial \Omega^+=\Sigma\cup \Sigma_0^+$, then $[\Sigma]=[\Sigma_0^+]=[S]\in H_2(M)$. On the other hand, for any surface $\Sigma$ in $[S]$ with $\Sigma\cap \Sigma_0=\emptyset$, the converse is also true, i.e. $\Sigma\subset \wh{M}$ separates $\Sigma_0^+$ and $\Sigma_0^-$.

By using this simple observation, we will prove the following lemma.

\begin{lem} \label{nogaplem} [No Gap in  the CMC Spectrum] Let $M$ be a closed Riemannian $3$-manifold. Let $S$ be an embedded closed surface in $M$ with $[S]\neq 0\in H_2(M)$. Let $\Sigma_0$ be the area minimizing surface in $\Sigma_0$. If there exists an embedded $H_0$-surface $\Sigma_{H_0}$ in $[S]$ with $\Sigma_0\cap\Sigma_{H_0}=\emptyset$, then for any $H\in [0,H_0]$, there exists an embedded $H$-surface $\Sigma_H$ in $[S]$.
\end{lem}

\begin{pf} Let $\wh{M}$ be the metric completion of $M-\Sigma_0$ as discussed above. Then, $\partial \wh{M}=\Sigma_0^+\cup\Sigma_0^-$. Assuming the mean curvature vector $\mathbf{H}$ on $\Sigma_H$ points towards $\Sigma_0^+$, let $\Omega_{H_0}$ be the region separated by $\Sigma_0^+$ and $\Sigma_{H_0}$. If $\mathbf{H}$ points other direction, same argument works with the region separated by  $\Sigma_0^-$ and $\Sigma_{H_0}$.
	
Now, consider $\Omega_{H_0}$. Let $\Sigma$ be a surface in $\Omega_{H_0}$ separating $\Sigma_0^+$ and $\Sigma_{H_0}$. Let $\Omega$ be the region in $\Omega_{H_0}$ separated by $\Sigma$ and $\Sigma_{H_0}$. For $H\in [0,H_0]$, define a functional $\A_H=|\Sigma|+2H||\Omega||$ where $|.|$ represents area, and $||.||$ represents the volume.

\begin{figure}[h]
	
	\relabelbox  {\epsfxsize=5in
		
		\centerline{\epsfbox{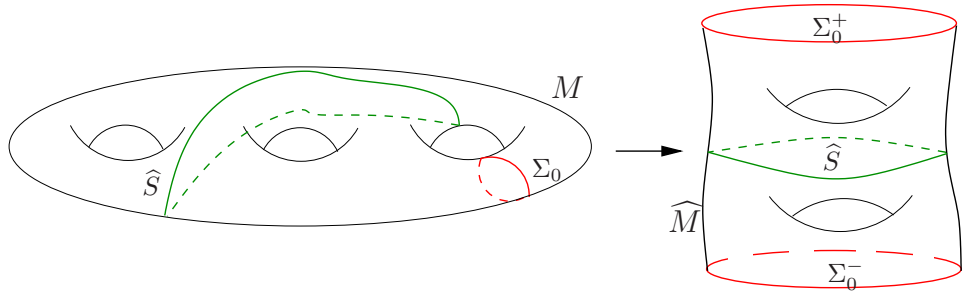}}}

	\relabel{1}{\footnotesize $\Sigma_0$} 
	\relabel{2}{\footnotesize $\wh{S}$} 
	\relabel{3}{$M$}
	\relabel{4}{$\wh{M}$} 
	\relabel{5}{\footnotesize $\Sigma_0^+$} 
	\relabel{6}{\footnotesize $\Sigma_0^-$}
	\relabel{7}{\footnotesize $\wh{S}$}

	\endrelabelbox
	
	\caption{\label{fig2} \small In the figure left, $\Sigma_0$ represents an embedded geodesic curve with nontrivial homology. $\wh{S}$ is a homologous curve to $\Sigma_0$.  In the figure right, we get a new manifold $\wh{M}$ after cutting $M$ along $\Sigma_0$. Any curve separating $\Sigma_0^+$ and $\Sigma_0^-$ in $\wh{M}$ is in the same homology class with $\Sigma_0$ in $M$.}
	
\end{figure}

By the standard theorems of geometric measure theory, for $H\in [0,H_0]$, $\A_H$ has a minimizer $\Sigma_H$ in $\Omega_{H_0}$. Furthermore, by maximum principle, $\Sigma_H$ is disjoint from $\partial \Omega_{H_0}=\Sigma_0^+\cup \Sigma_{H_0}$. Therefore, $\Sigma_H$ is a smoothly embedded $H$-surface in $\Omega_{H_0}$ by the regularity of the minimizers. Furthermore, any such $\Sigma_H$ is separating in $\wh{M}$ by construction, therefore $\Sigma_H\in [S]$. Hence, for $H\in [0,H_0]$  there exists an embedded $H$-surface $\Sigma_H$ in $[S]$. The proof follows.
\end{pf}

\begin{rmk} Notice that the lemma above implies that if there exists an $H_0$-surface in a homology class $[S]$, then the whole interval $[0,H_0]$ is achieved. In particular, \textit{there is no gap in the CMC spectrum of the homology class}. 	
\end{rmk}

\subsection{A Lower Bound for $\wh{H}_{[S]}$} \

Let $\wh{H}_{[S]}=\sup\{|H| \mid \exists \ \Sigma_H \in [S]\}$. After the lemma above, we only need to establish a good lower bound for $\wh{H}$ to obtain our main result. In order to establish such lower bound, we need to study the isoperimetric profile of the homology class.

Let $M, S, \Sigma_0,\wh{M}$ be described as above. We will define the \textit{isoperimetric profile function} $\I_{[S]}$ of the homology class $[S]$ as follows.  Let $\T$ be an embedded surface in $\wh{M}$ separating $\Sigma_0^+$ and $\Sigma_0^-$. Let $\Delta$ be the region in $\wh{M}$ separated by $\Sigma_0^+$ and $\Sigma$. Let $K\in [0,||M||]$. Then define 
$$\I_{[S]}(K)=inf\{|\T| \mid ||\Delta||=K \}$$

Hence, we defined a function $\I_{[S]}: [0, ||M||]\to \BR^+$. This function will be key to get our lower bound. First, we will show that this function is continuous.

\begin{lem} \label{ctslem} The isoperimetric profile function $\I_{[S]}(K)$ is continuous.
\end{lem}

\begin{pf} By \cite{Ro}, we have the continuity of the isoperimetric function for closed manifolds. We will use this result for our case as follows. Let $D_1$ be the diameter of $\Sigma_0$. Let $D_2$ be the distance between $\Sigma_0^+$ and $\Sigma_0^-$ in $\wh{M}$. Let $n>>\frac{D_1}{D_2}$. Now, let $\wh{M}^n$ be the $n$-cover of $\wh{M}$. In particular, we glue $n$ copies of $\wh{M}$ to each other one by one as follows. Let $\wh{M}_i$ be a copy of $\wh{M}$ with boundary components $\Sigma_i^\pm$. Then, glue the lower end $\Sigma_i^-$ of one to the upper end $\Sigma_{i+1}^+$ of the next. Hence, $\wh{M}^n$ contains $n$ copies of $\wh{M}$, and $2$ boundary components $\Sigma_1^+$ and $\Sigma_n^-$. Now, cap off $\wh{M}^n$ by gluing smoothly small diameter handlebodies $N^+$ and $N^-$ to $\partial \wh{M}^n$. We obtain a smooth closed manifold $\wt{M}$ which can be decomposed as $n$ copies of $\wh{M}$ and 2 handlebodies $N^+$ and $N^-$.

As $\wt{M}$ is a closed $3$-manifold, the isoperimetric profile function $\I_{\wt{M}}$ of $\wt{M}$ is continuous by \cite{Ro}. In particular, $\I_{\wt{M}}: [0,\|\wt{M}\|]\to \BR$. As we choose $n$ sufficiently large, in the middle segment of the image of $\I_{\wt{M}}$ will be periodic. In particular, let  $K_j=\|N^+\|+j\|M\|$. Then for any $1< j<j'< n$, we will have  $\I_{\wt{M}}(K_j)=\I_{\wt{M}}(K_{j'})=|\Sigma_0|$ by construction. On the other hand, one of these periodic segments will exactly match the image of $\I_{[S]}$ for $\wh{M}$ by construction. In particular, when $j\in\BN$ sufficiently large, say $\sim \frac{n}{2}$, we will have $\I_{\wt{M}}(\|N^+\|+j\|M\|+K)=\I_{[S]}(K)$. This shows that $\I_{[S]}$ is a continuous function. The proof follows.
\end{pf}

As mentioned earlier, in \cite{Ro}, there are several results given for the isoperimetric profile of closed $3$-manifold $M$. In the proof of the lemma above, we show that the graph of $\I_{[S]}$ is the same with a portion of graph of the isoperimetric profile of the closed $3$-manifold $\wt{M}$. Hence, we can adapt these results to our setting $\wh{M}$. 

In particular, for any $K\in [0,||M||]$, the right and left derivatives of the function $\I_{[S]}(K)$ exists, i.e. $\I_{[S]}^{'+}(K)$ (right derivative) and $\I_{[S]}^{'-}(K)$ (left derivative). Furthermore, there exists an isoperimetric surface $\T_K^\pm$ with mean curvature $H^\pm=\frac{\I_{[S]}^{'\pm}(K)}{2}$. By using this, we will get our lower bound. Note that the derivative can be negative, which means the mean curvature vector changed direction.

Now, we will use a simple observation on the graph of $\I_{[S]}$ to prove the following lemma. One can think of the following lemma as a version of the mean value theorem. Recall that girth of a homology class is defined as $g([S])=\max_K \I_{[S]}(K)$.

\begin{lem} \label{lowerlem} Let $C_S=\dfrac{g([S])-|\Sigma_0|}{||M||}$. Then, there exists  $ K_1\in (0,||M||)$ such that $\max\{|\I_{[S]}^+(K_1)|,|\I_{[S]}^-(K_1)|\}\geq 2C_S$.  
\end{lem}

\begin{pf} By Lemma \ref{ctslem}, the isoperimetric profile function $\I_{[S]}(K)$ is continuous. As $\max \I_{[S]}=g([S])$, say $\I_{[S]}(K_0)=g([S])$. Recall that by construction, $\min \I_{[S]} = \I_{[S]}(0)=\I_{[S]}(||M||)=|\Sigma_0|$. 
	
As $\I_{[S]}$ is continuous, there must be a point $K_1\in (0, K_0)$ with (right or left) derivative at least $\dfrac{\I_{[S]}(K_1)-\I_{[S]}(0)}{K_1-0}$. Similarly, there must be a point $K_1'\in (K_0,||M||)$ with (right or left) derivative at most $\dfrac{I_{[S]}(||M||)-I_{[S]}(K_0)}{||M||-K_0}$. 

As $K_0\in (0,||M||)$, $\max\{K_0, (||M||-K_0)\}\geq \dfrac{||M||}{2}$.

Then, $\max\{\I_{[S]}^{'\pm}(K_1), \I_{[S]}^{'\pm}(K_1')\}\geq \dfrac{2(g([S])-|\Sigma_0|)}{||M||}=2C_S$. 

The proof follows.
\end{pf}

Hence, the existence of such a left or right derivative in the graph of $\I_{\wh{M}}$ gives the following corollary.

\begin{cor} \label{lowercor} Let $\wh{H}=\max\{|H| \mid \exists \Sigma_H \in [S]\}$. Then,  $\wh{H}\geq C_{[S]}$.
\end{cor}

\begin{pf} Notice that by Lemma \ref{nogaplem}, we only need to show that there exists an $H$-surface $\Sigma_H$ in $[S]$ where $H\geq C_{[S]}\geq 0$.	
	
By Lemma \ref{lowerlem}, there exists $K_1\in[0,||M||]$ such that $|\I_{S}^{'\pm}(K_1)|\geq 2C_S$. Recall that by \cite{Ro}, there is an isoperimetric surface $T^\pm_{K_1}$ such that the mean curvature along $T^\pm_{K_1}$ is $H^\pm_{K_1}= \dfrac{\I_{S}^{'\pm}(K_1)}{2}$. So, either $T^+_{K_1}$ or $T^-_{K_1}$ is an $H$-surface in $\wh{M}$ with $H\geq C_{[S]}$. By construction, $T^\pm_{K_1}$ is in the homology class of $[S]$. The proof follows.
\end{pf}

\begin{rmk} [Orientation] Note that if the derivative $\I_{[S]}^{'\pm}(K_1)$ is positive, then depending on the orientation we choose, the mean curvature vector along $T^\pm_{K_1}$ points toward $\Sigma_0^+$ or $\Sigma_0^-$. Hence, if the derivative is negative, the mean curvature vector points the opposite surface. Recall that if we use negative mean curvature surface $\T$, we use the region between $\T$ and $\Sigma_0^-$ to define the functional $\A_H$ in Lemma \ref{nogaplem}. 
\end{rmk}

\begin{rmk} [girth vs. width]  It is clear that for a given manifold $M$, and a given homology class $[S]\in H_2(M)$, it is not easy to compute the girth $g([S])$ in general. However, we have an upper bound for $g([S])$ which is relatively easy to compute: The width $w([S])$ of the homology class, which is the area of the min-max surface in $[S]$. See Section \ref{minmaxsec} for further details.	
\end{rmk}

Now, we can state our main result.

\begin{thm} \label{mainthm}	Let $M$ be a closed Riemannian $3$-manifold. Let $S$ be a closed embedded surface in $M$, where $S$ has nontrivial homology. Then, there exists $\widehat{H}\geq C_{S}\geq 0$ such that for any $H\in [0,\wh{H}]$, there exists a smoothly embedded $H$-surface $\Sigma_H$ in the homology class of $S$.	
\end{thm}

\begin{pf} By Corollary \ref{lowercor}, we have an isoperimetric surface $\T_{K_1}$ in $[S]$ with constant mean curvature $H\geq C_{[S]}$. Then, by Lemma \ref{nogaplem}, there exists $\widehat{H}\geq C_{[S]}\geq 0$ such that for any $H\in [0,\wh{H}]$, there exists a smoothly embedded $H$-surface $\Sigma_H$ in $[S]$. The proof follows.
\end{pf}
	
%The proof follows from Lemma \ref{nogaplem}, and Corollary \ref{lowercor}.

%By \cite{Ro}, for any $K\in [0,\|M\|]$, we will have an isoperimetric surface $\T_K^+$ and $\T_K^-$ with mean curvature $\I_{[S]}^{'+}(K)$  and $\I_{[S]}^{'-}(K)$ respectively. Then,

\section{Concluding Remarks}

\subsection{Minimizing $H$-Surfaces vs. Isoperimetric Surfaces} \

\vspace{.2cm}

So far, we used two different minimization techniques to construct CMC Surfaces i.e. minimizing $H$-surfaces and Isoperimetric Surfaces. In particular, if we fix $H$, and minimize $\A_H$, we get minimizing $H$-surface $\Sigma_H$. If we fix volume $||\Omega||=K$ and minimize area $|\T|$, we get isoperimetric surface $\T_K$. Hence, the natural question is the following:

\begin{question} What is the relation between minimizing $H$-surfaces $\Sigma_H$ and Isoperimetric Surfaces $\T_K$. Is there a correspondence between $\Sigma_H$ and $\T_K$?
\end{question}

The following lemma is a starting point to address this question.

\begin{lem} Any $H$-minimizing surface is an isoperimetric surface. In particular, $\Sigma_H=\T_K$ for $K=||\Omega_H||.$
\end{lem}

\begin{pf} Fix $H_0>0$. Let $\Sigma_{H_0}$ be the minimizing $H_0$-surface in $\wh{M}$. We claim that $\Sigma_{H_0}$ is isoperimetric surface for $K_0=||\Omega_{H_0}||$ where $\Omega_{H_0}$ is the region between $\Sigma_0^+$ and $\Sigma_{H_0}$. 
	
In particular, consider the functional $\wh{\A}_{H_0}(\Sigma)=|\Sigma|-2H_0|\Omega|$ (See Section \ref{minHsec}). Consider all surfaces $\Sigma$ in $\wh{M}$ with $||\Omega||=K_0$ where $\Omega$ is the region between $\Sigma$ and $\Sigma_0^+$. Recall that $\Sigma_{H_0}$ minimizes $\wh{\A}_{H_0}$ among all such surfaces. Then, $\wh{\A}_{H_0}(\Sigma_{H_0})\leq \wh{\A}_{H_0}(\Sigma)$ for any $\Sigma \in [S]$.
	
Let $\Sigma'$ be such that $||\Omega'||=K_0$. 	Then, we have 
	
	$\wh{\A}_{H_0}(\Sigma_{H_0})=|\Sigma_{H_0}|-2H_0 K_0$ while $\wh{\A}_{H_0}(\Sigma')=|\Sigma'|-2H_0 K_0$. 
	
	Then, we have $$|\Sigma_{H_0}|-2H_0K_0\leq |\Sigma|-2H_0K_0$$ which implies $|\Sigma_{H_0}|\leq |\Sigma'|$. The proof follows. 
\end{pf}

Then, the natural followup question is:

\begin{question} Is every isoperimetric surface $\T_K$ in $\wh{M}$ is a minimizing $H$-surface in $[S]$ for $H=H(\T_K)$?
\end{question}

We will have a partial answer to this question. Now, consider the correspondence between the mean curvature $H$, and the volume $K$ for isoperimetric surfaces and minimizing $H$-surfaces. Define the following "multivalued" functions $\h(K)$ and $\kk(H)$: $$\h(K)=H(\T_K^\pm) \quad \mbox{where} \quad K\in [0,\|M\|]$$

where $T^\pm_K$ is the isoperimetric surface in $\wh{M}$ with $||\Omega||=K$. Here, $\T^+_K$ and $\T^-_K$ might be two different isoperimetric surfaces with $\|\Omega^+_K\|=\|\Omega^-_K\|=K$. $H(T^+_K)$ is the constant mean curvature of $\T^+_K$, and $H(\T^-_K)$ is defined similarly. In particular, $H(\T_K^\pm)=\dfrac{\I_{[S]}^{'\pm}(K)}{2}$ where $\I_{[S]}$ is the isoperimetric profile of the homology class $[S]$. Notice that whenever $\I_{[S]}$ is not smooth, we would have 2 different isoperimetric surface. In particular, if $\I_{[S]}$ is not smooth at $K_0\in (0,\|M\|)$, then we would have two different isoperimetric surfaces $T_{K_0}^+$ and $T_{K_0}^-$  with constant mean curvatures $\frac{\I_{[S]}^{'+}(K_0)}{2}$ and $\frac{\I_{[S]}^{'-}(K_0)}{2}$ respectively.  In particular, by \cite{Ro} and Lemma \ref*{ctslem}, $\h(K)=\frac{\I_{[S]}'(K)}{2}$ when $\I_{[S]}'(K)$ is defined. Hence, $\h(K)$ is multivalued when $\I_{[S]}'(K)$ is not defined, i.e. $\I_{[S]}$ is not smooth at $K\in (0,\|M\|)$. 

%In other words, whenever the isoperimetric profile function $\I_{[S]}$ is not smooth, then the function $\h(K)$ is multivalued.

Now, define the "inverse" multivalued function 

$$\kk(H)=||\Omega_H|| \quad \mbox{where} \quad H\in[0,\wh{H}_{[S]}]$$

where $\Omega_H$ is the region between minimizing $H$-surface $\Sigma_H$ and $\Sigma_0^+$. Recall that every minimizing $H$-surface is an isoperimetric surface for the corresponding volume by Lemma \ref{nogaplem}. Notice that for the same $H_0\in[0,\wh{H}_{[S]}]$ there might be more than one minimizing $H_0$-surface as they are defined as the minimizers of the functional $\A_{H_0}$. Assume that $\A_{H_0}$ have more than one minimizers $\Sigma_{H_0}^+$ and $\Sigma_{H_0}^-$. Then by Lemma \ref{uniqthm}, we have $\|\Omega_{H_0}^+\|\neq \|\Omega_{H_0}^-\|$, which implies that the function $\kk(H)$ is multivalued at $H_0$. 

Now, by using the result in the appendix, we can conclude the following.

\begin{cor} For $H\in [0,\wh{H}_{[S]}]$, $\kk(H)$ is increasing. 		
\end{cor}

\begin{pf} By Theorem \ref{uniqthm}, when $H_1<H_2$, then $\Omega_{H_1}\subset \Omega_{H_2}$. This would imply $\|\Omega_{H_1}\|<\|\Omega_{H_2}\|$. Therefore, we have $\kk(H_1)<\kk(H_2)$. Notice that as $\kk(H)$ is multivalued function, there is an ambiguity here. However, this does not effect our result. In particular, if $\kk(H)$ is multivalued at $H_0$, and $H_1<H_0<H_2$, then by Theorem \ref{uniqthm}, we have $\Omega_{H_1}\subset \Omega_{H_0}^-\subset \Omega_{H_0}^+\subset \Omega_{H_2}$ (See Figure \ref{fig3}-right). Hence, we have $\kk(H_1)<\kk(H_0^-)<\kk(H_0^+)<\kk(H_2)$. The proof follows.
\end{pf}

%What is the motivation for this question/conjecture? First notice that, we are only considering $[0,\wh{H}_{[S]}]$ interval for the function $\kk(H)$. This is intuitively corresponding to the first "quarter" of the isoperimetric profile $\I_{[S]}$. This is a restricted and somewhat special region. Recall Lemma \ref{nogaplem}. In the proof of that lemma, if we fix our barrier surface as $\Sigma_{\wh{H}_{[S]}}$, we would get minimizing $H$-surfaces in the region $\Omega_{\wh{H}_{[S]}}$ which is the region between $\Sigma_0^+$ and $\Sigma_{\wh{H}_{[S]}}$ for any $H\in[0,\wh{H}_{[S]}]$. Furthermore, in the same lemma, we proved that these embedded minimizing $H$-Surfaces are pairwise disjoint for different $H$, i.e. $\Sigma_{H_1}\cap\Sigma_{H_2}=\emptyset$. This observation suggest that when $H$ increase from $0$ to $\wh{H}_{[S]}$, then the volume of the region $\|\Omega\|$ increases as well. In other words, there is an ordered structure among these minimizing $H$-surfaces for $H\in[0,\wh{H}_{[S]}]$. If this is true, then it would imply that the function $\kk(H)$ is non-decreasing?

%GENERIC UNIQUENESS OF $\Sigma_H$ FOR $[0,\wh{H}_{[S]}]$ WHAT DOES THIS SAY ABOUT $\kk(H)$.

Notice that $\kk(H)$ and $\h(K)$ functions are somewhat "inverse" functions. By using this relation, what can we say about $\h(K)$? Or in particular,

\begin{question} For a given embedded surface $S\subset M$ with nontrivial homology, what can we say about the isoperimetric profile function $\I_{[S]}(K)$, or its derivative $\I_{[S]}^{'\pm}(K)=2\h(K)$?	
\end{question}

%Regularity of $\I_{[S]}$? If nowhere smooth, all these ideas collapse..

%{\large By using these or different ideas, can we improve the lower bound for $\wh{H}$?}

\subsection{Min-Max Surfaces and Improving the Lower Bound $C_{[S]}$} \label{minmaxsec} \

\vspace{.2cm}

In this paper, our methods produce locally minimizing CMC surfaces. In particular, all the surfaces $\Sigma_H$ we produced in this paper are {\em stable CMC surfaces}. However, there are other ways to obtain CMC surfaces like min-max method. Min-max method generates {\em unstable CMC surfaces} which are different than what we constructed so far. Recently, Zhou and Zhu gave a remarkable result by showing the existence of almost embedded $H$-surfaces in any closed $3$-manifold for any constant $H$ \cite{ZZ}. By adapting their techniques in our setting for $\wh{M}$, it might be possible to improve our lower bound for $\wh{H}_{[S]}$. In particular, we pose the following conjecture, which is basically a generalization of our main result:

\begin{conj} Let $M$ be a closed $3$-manifold. Let $S$ be an embedded surface in $M$ where $S$ has nontrivial homology. Then, there exists $\widehat{H}_{[S]}\geq \wh{C}_{[S]}\geq 0$ such that for any $H\in [0,\wh{H}]$, there exists a smoothly embedded $H$-surface $\Sigma_H$ in the homology class of $S$.	
\end{conj}

Here, $\wh{C}_{[S]}=\dfrac{w([S])-|\Sigma_0|}{||M||}$ where $|.|$ represents area, and $||.||$ represents volume. $w([S])$ is {\em the width of the homology class $[S]$}, which is the area of the min-max surface $\wt{\Sigma}$ in the homology class $[S]$, then $w([S])=|\wt{\Sigma}|$. 

In particular, let $\wh{M}$ be defined as before. Let $$X_{[S]}= \{ f:\wh{M}\to [0,1] \mid f \mbox{ is a Morse function where } f(\Sigma_0^+)=1 \mbox{ and } f(\Sigma_0^-)=0 \}$$
Then, for any $f\in X_{[S]}$, define sweepout $\Lambda_f=\{\Sigma^f_t=f^{-1}(t)\}$. Then, define the min-max surface $\wt{\Sigma}$ of the homology class $[S]$ as the minimum of the maximal slices of these sweepouts, i.e. $|\wt{\Sigma}|=\min_{f\in X_{[S]}}\max_t |\Sigma^f_t|$.

Note that it is not hard to show the existence of such min-max surface $\wt{\Sigma}$ in $\wh{M}$ as defined above by using the technique in the proof of Lemma \ref{ctslem}. In particular, consider the closed $3$-manifold $\wt{M}$ defined in the proof of Lemma \ref{ctslem}. We have a min-max surface $\wt{\Sigma}$ in the closed $3$-manifold $\wt{M}$ \cite{ZZ}. Then, by the construction of $\wt{M}$, it is not hard to see that the same surface $\wt{\Sigma}$ in the closed manifold $\wt{M}$ corresponds to the min-max surface $\wt{\Sigma}$ in $\wh{M}$ as defined above.

%PROOF OF EXISTENCE OF MIN MAX SURFACE. wt{M}= MX [0, K] U M+ U M- extend sweepouts. [ZZ] IMPLIES EXISTENCE.

This conjecture is a generalization of our result (Theorem \ref{mainthm}) because of the following lemma, i.e. $\wh{C}_{[S]}\geq C_{[S]}$.

\begin{lem} Let $M$ be a closed $3$-manifold. Let $S$ be an embedded surface in $M$ where $S$ has nontrivial homology. Then, $w([S])\geq g([S])$.
\end{lem}

\begin{pf} 	Let $\I_{[S]}(K_0)=max \I_{[S]}(K)$. Then, for any sweepout $\Lambda_f$ consider the surface $S^f_{K_0}$ where $||\Omega^f_{K_0}||=K_0$. Then, $|S^f_{K_0}|\geq \I_{[S]}(K_0)$ by the definition of $\I_{[S]}(K)$. Since maximal slice $S^f$ in $\Lambda_f$ has greater area than $S^f_{K_0}$, we have $\max\Lambda_f\geq \max \I_{[S]}(K)$ for any sweepout $\Lambda_f$. This implies $w([S])=\min\max_{\Lambda_f}|S|=\min (|S^f|)\geq \max \I_{[S]}(K_0)=g([S])$	
\end{pf}

Here, we also conjecture that $\wh{C}_{[S]}$ is the best lower bound for $\wh{H}_{[S]}$ in general. However, this is a much harder problem than the generalization above, as it needs a much deeper understanding of the original question.

\section{Appendix}

In this part, we give an interesting observation related to relation between isoperimetric surfaces and minimizing $H$-surfaces. One can see this result as an extension of the Lemma \ref{nogaplem}, {\em No Gap Lemma}. For the brevity of the discussion, we postponed this result to this section.

\begin{thm} \label{uniqthm} Let $M$ be a closed Riemannian $3$-manifold. Let $S$ be an embedded closed surface in $M$ with $[S]\neq 0\in H_2(M)$. Let $\Sigma_0$ be the area minimizing surface in $\Sigma_0$. If there exists an embedded $H_0$-surface $\Sigma_{H_0}$ in $[S]$ with $\Sigma_0\cap\Sigma_{H_0}=\emptyset$, then 

\begin{itemize}
	\item For any $H\in [0,H_0]$, there exists an embedded $H$-surface $\Sigma_H$ in $[S]$.
	
	\item For any $H_1\neq H_2 \in [0,H_0]$, we have $\Sigma_{H_1}\cap\Sigma_{H_2}=\emptyset$.
	
	\item For a generic $H\in [0,H_0]$, there exists a unique minimizing $H$-Surface $\Sigma_H$ in $\Omega_{H_0}$.	
	
\end{itemize}
		
\end{thm}

\begin{pf} The first part of the theorem was already proven in Lemma \ref{nogaplem}. We will prove the following two parts in three steps.
	
\vspace{.2cm}
	
\noindent {\bf Step 1:} Let $H_1\neq H_2 \in [0,H_0]$. Then $\Sigma_{H_1}\cap\Sigma_{H_2}=\emptyset$

\vspace{.2cm}

\noindent {\em Proof of Step 1:} Let $\Sigma_{H_1}$ and $\Sigma_{H_2}$ be minimizers of $\wh{\A}_{H_1}$ and $\wh{\A}_{H_2}$ in $\Omega_{H_0}$ respectively (See Section \ref{minHsec}). Let $H_1<H_2$. We claim that $\Sigma_{H_1}\cap \Sigma_{H_2}=\emptyset$. 

For $i=1,2$, let $\Omega_{H_i}$ be the closed region between $\Sigma_{H_i}$ and $\Sigma_0^+$. We will show that $\overline{\Omega}_{H_1}\subset int(\Omega_{H_2})$.

Assume that $\Sigma_{H_1}\cap \Sigma_{H_2}\neq \emptyset$. Let $\Delta=\Omega_{H_1}\cap\Omega_{H_2}$. Let $S_1=\Sigma_{H_1}\cap \Omega_{H_2}$ and $T_1=\Sigma_{H_1}- int(\Omega_{H_2})$. Similarly, let $S_2=\Sigma_{H_2}- int(\Omega_{H_2})$ and $T_1=\Sigma_{H_1}\cap \Omega_{H_2}$. Let $X_1=\Omega_{H_1}-\Omega_{H_2}$, and $X_2=\Omega_{H_2}-\Omega_{H_1}$. See Figure \ref{fig3}-left.

\begin{figure}[h]
	\begin{center}
		$\begin{array}{c@{\hspace{.4in}}c}

		\relabelbox  {\epsfxsize=1.7in \epsfbox{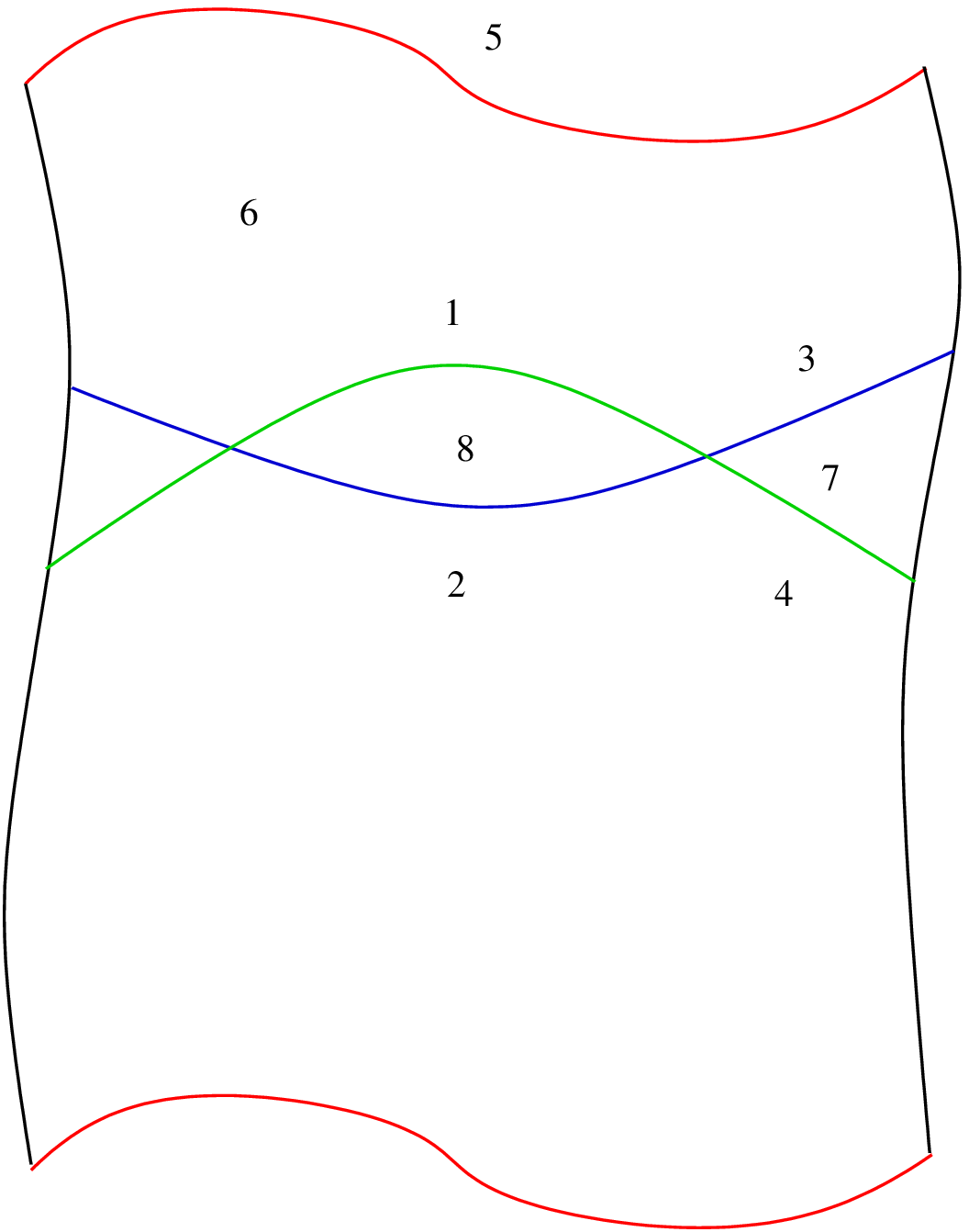}}
		\relabel{1}{\footnotesize $T_2$} 
		\relabel{2}{\footnotesize $T_1$} 
		\relabel{3}{$S_1$}
		\relabel{4}{$S_2$} 
		\relabel{5}{\footnotesize $\Sigma_0^+$} 
		\relabel{6}{\footnotesize $\Delta$}
		\relabel{7}{\footnotesize $X_2$}
		\relabel{8}{\footnotesize $X_1$}	   \endrelabelbox &
		
		\relabelbox  {\epsfxsize=1.7in \epsfbox{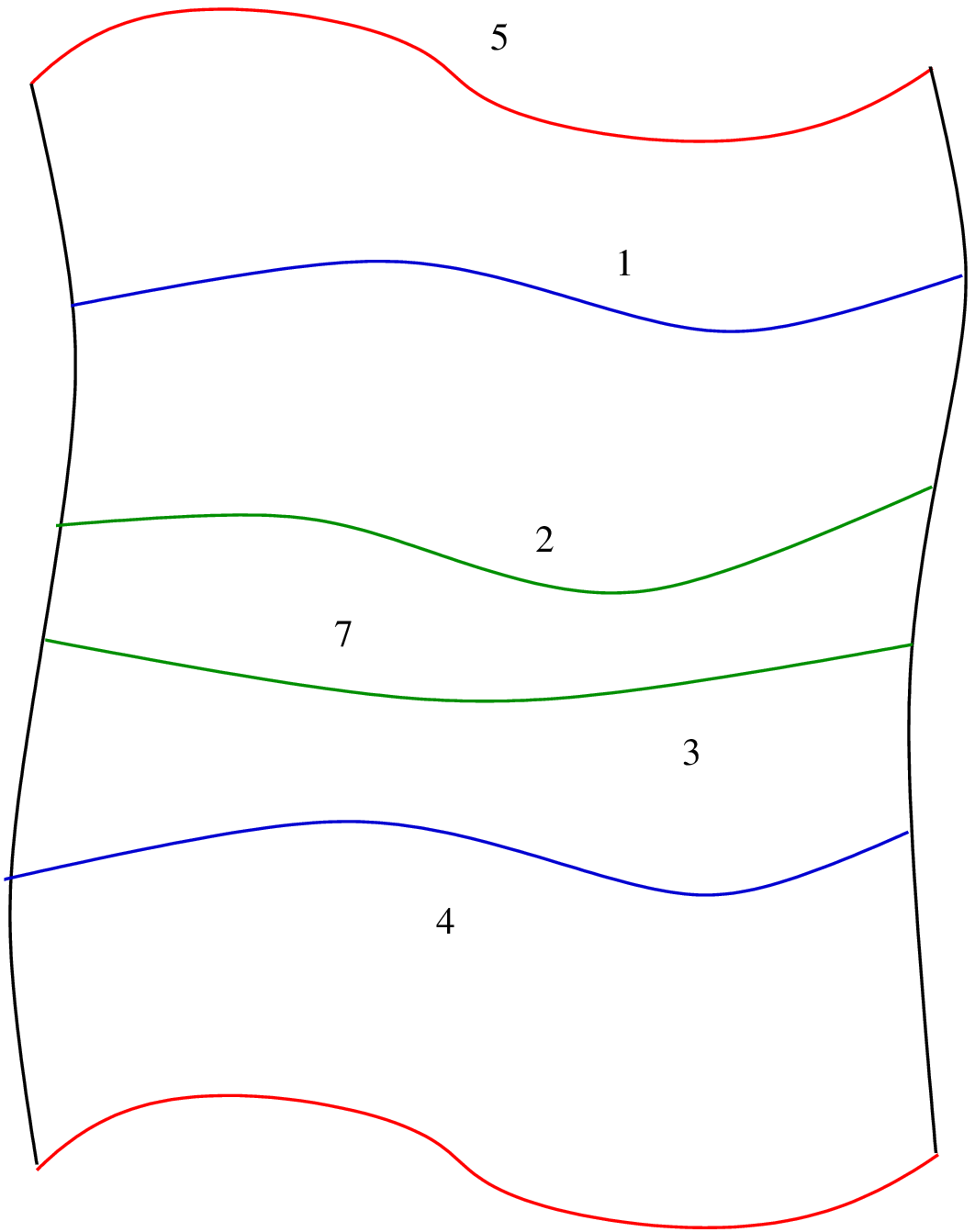}}
		\relabel{1}{\footnotesize $\Sigma_{H_1}$} 
		\relabel{2}{\footnotesize $\Sigma_{H_0}^-$}  
		\relabel{3}{\footnotesize $\Sigma_{H_0}^+$} 
		\relabel{4}{\footnotesize $\Sigma_{H_2}$}
		\relabel{5}{\footnotesize $\Sigma_0^+$}
		%\relabel{6}{\footnotesize $\Sigma_0^-$}
		\relabel{7}{\footnotesize $\Delta_{H_0}$}	\endrelabelbox \\
		
	\end{array}$
	
\end{center}

\caption{ \label{fig3} \footnotesize In the figure left, $\Sigma_{H_i}=S_i\cup T_i$ and $\Omega_{H_i}=X_i\cup \Delta$. In the figure right, for $H_1<H_0<H_2$, we have $\Sigma_{H_1}>\Sigma_{H_0}^->\Sigma_{H_0}^+>\Sigma_{H_2}$. $\Delta_{H_0}$ is the canonical region for $H_0$ with $\partial \Delta_{H_0}=\Sigma_{H_0}^+\cup\Sigma_{H_0}^-$.}
\end{figure}

Now, let $\Sigma_1'=S_1\cup T_2$ and $\Sigma_2'=S_2\cup T_1$. As $\Sigma_{H_i}$ is minimizer for $\wh{\A}_{H_i}$, we have the following: 
$$\wh{\A}_{H_i}(\Sigma_{H_i})\leq \wh{\A}_{H_i}(\Sigma_i') \quad \mbox{for} \quad i=1,2$$ 

This implies

\vspace{.2cm}

$\wh{\A}_{H_1}(\Sigma_{H_1})\leq \wh{\A}_{H_1}(\Sigma_1') \quad \Rightarrow \quad 
|\Sigma_{H_1}|-2H_1||\Omega_{H_1}|| \leq |\Sigma_1'|-2H_1||\Omega_1'||$

\vspace{.2cm}

Then, we have 

\vspace{.2cm}

$(|S_1|+ |T_1|) - 2H_1 (\|\Delta\|+\|X_1\|) \leq (|S_1|+|T_2|) - 2H_1(\|\Delta\|)$
$$|T_1| -2H_1\|X_1\|\leq |T_2|  \quad \quad (\star)$$

On the other hand, 

\vspace{.2cm}

$\A_{H_2}(\Sigma_{H_2})\leq \A_{H_2}(\Sigma_2') \quad \Rightarrow \quad 
|\Sigma_{H_2}|-2H_2||\Omega_{H_2}|| \leq |\Sigma_2'|-2H_2||\Omega_2'||$

\vspace{.2cm}

Then, we have 

\vspace{.2cm}

\noindent $(|S_2|+ |T_2|) - 2H_2 (\|\Delta\|+\|X_2\|) \leq (|S_2|+|T_1|) + 2H_2(\|\Delta\|+\|X_1\|+\|X_2\|)$
$$|T_2||  \leq |T_1| -2H_2\|X_1\| \quad \quad (\star \star)$$

By combining $(\star)$ and $(\star \star)$, we have
$$2H_2||X_1||\leq |T_1|-|T_2|\leq 2H_1||X_1||$$

As we assumed $||X_1||> 0$, this gives $H_2\leq H_1$. However, our assumption was $0<H_1<H_2$. This is a contradiction. Step 1 follows. \hfill $\Box$

\vspace{.2cm}

\noindent {\bf Step 2:} For any $H\in [0,H_0]$, either there exists a unique minimizing $H$-surface $\Sigma_H$, or there are two canonical minimizing $H$-surfaces $\Sigma_H^+$ and $\Sigma_H^-$ in $\wh{M}$.

\vspace{.2cm}

\noindent {\em Proof of Step 2:} Basically, we will adapt the ideas in \cite{Co2} to this setting. By Lemma \ref{nogaplem}, for $H_1<H_2$, $\Sigma_{H_1}\cap\Sigma_{H_2}=\emptyset$ and $\Omega_{H_1}\subset \Omega_{H_2}$. Now, for any $H_0\in (0,\wh{H}_{[S]})$, define $\Sigma_{H_0}^+$ and $\Sigma_{H_0}^-$ as follows. Take a monotone decreasing sequence $H_n\searrow H_0$. Consider the sequence $\{\Sigma_{H_n}\}$. Then, by compactness theorem of geometric measure theory, we will get a minimizing $H_0$-surface in the limit, i.e. $\lim \Sigma_{H_n}=\Sigma_{H_0}^+$. 

We claim that $\Sigma_{H_0}^+$ is canonical, i.e. independent of the sequence $\{\Sigma_{H_n}\}$. In order to see this, take another monotone sequence $H_m\searrow H_0$, and corresponding sequence $\{\Sigma_{H_m}\}$. We claim that $\lim \Sigma_{H_m}=\Sigma_{H_0}^+$ again. Assume not. Let $\lim \Sigma_{H_m}=\Sigma'$ which is another $H_0$-minimizing surface. By Lemma \ref{nogaplem}, we have an ordered structure $\Omega_1\supset\Omega_2\supset...\supset\Omega_m\supset...$. Let $\Omega_{H_0}^+$ be the region between $\Sigma_0^+$ and $\Sigma_{H_0}^+$, and let $\Omega'$ be the region between $\Sigma_0^+$ and $\Sigma'$. We claim that $\Omega_{H_0}^+=\Omega'$. Assume not. If $\Omega'$ is not in $\Omega_{H_0}^+$, this would imply for sufficiently large $n_0$, $\Sigma_{H_{n_0}}\cap \Sigma'\neq \emptyset$. However, this contradict to Step 1. Similarly, if $\Omega_{H_0}^+$ is not in $\Omega'$, this would imply for sufficiently large $m_0$, $\Sigma_{H_{m_0}}\cap \Sigma_{H_0}^+\neq \emptyset$. However, this contradict to Step 1 again. Similarly, we can define $\Sigma_{H_0}^-$ as $\lim \Sigma_{H_n}$ where $H_n\nearrow H_0$ a mononotone increasing sequence. See Figure \ref{fig3}-right. Same reasoning shows that if $\Sigma_{H_0}^+=\Sigma_{H_0}^-$, then there exists a unique $H_0$-minimizing surface $\Sigma_{H_0}$. Step 2 follows. \hfill $\Box$

\vspace{.2cm}
	
\noindent {\bf Step 3:} For a generic $H\in [0,H_0]$, there exists a unique minimizing $H$-Surface $\Sigma_H$ in $\wh{M}$.

\vspace{.2cm}

\noindent {\em Proof of Step 3:} Now, for any $H\in (0,H_0)$, define a canonical neighborhood $\Delta_H=[\Sigma_H^-,\Sigma_H^+]$, i.e. the closed region between $\Sigma_H^-$ and $\Sigma_H^+$. In particular, if $\Sigma_H^-=\Sigma_H^+$, then $\Delta_H=\Sigma_H$. Take a finite transversal arc $\beta$ intersecting all $\{\Sigma_H^\pm\}$ for $H\in [0,H_0]$. Let $s_H$ be the length of $\beta\cap\Delta_H$. Notice that if $s_H=0$, then there exists a unique minimizing $H$-surface $\Sigma_H$ in $\Omega_{H_0}$.  Similarly, if $s_H>0$, we have more than one minimizing $H$-surfaces $\Sigma_H^+$ and $\Sigma_H^-$. We call $s_H$ the {\em thickness} of $\Delta_H$. 

Notice that $\sum s_H \leq\|\beta\|$ by construction. As the sum is finite, for only countably many $H\in [0,H_0]$, $s_H$ can be nonzero. This implies for all but countably many of $H\in [0,H_0]$, we have $s_H=0$, and hence, unique minimizing $H$-surface. Step 3, and the proof of the theorem follows.	
\end{pf}

%FIGURE: $H_1<H_2 \Rightarrow \Sigma_1>\Sigma_2$. Also, $H_i\searrow H_0 \Rightarrow \Sigma_i\to \Sigma_0^-$ and $H_i\nearrow H_0 \Rightarrow \Sigma_i\to \Sigma_0^+$. Then, $\Sigma_0^+\cap \Sigma_0^-=\emptyset$ and $\Sigma_0^->\Sigma_0^+$. Draw FIGURE!!!

\end{document}